# Integral points on the rational curve

$$y = \frac{x^2 + bx + c}{x + a}; \quad a,b,c \text{ integers.}$$

**Konstantine Zelator**

**Mathematics**

**University of Wisconsin - Marinette**

**750 W. Bayshore Street**

**Marinette, WI 54143-4253**

*Also:*

**Konstantine Zelator**

**P.O. Box 4280**

**Pittsburgh, PA 15203**

Email Addresses:   1)  konstantine_zelator@yahoo.com

2)  spaceman@pitt.edu

## 1. Introduction

In this work, we investigate the subject of integral points on rational curves of the form

$y = \dfrac{x^2 + bx + c}{x + a}$; where $a, b, c$ are fixed or given integers. An integral point is a point *(x,y)* with both coordinates *x* and *y* being integers. Not that if *($i_1$, $i_2$)* is an integral point with $i_1 \neq -a$; on the graph of the above rational function. Then *($i_1$, $i_2$)* also satisfies the algebraic equation, of degree *2*; $x^2 - xy + bx - ay + c = 0$. And conversely, if *($i_1$, $i_2$)* is an integer pair with $i_1 \neq -a$; satisfying this algebraic equation. It then lies on the graph of the rational function $y = \dfrac{x^2 + bx + c}{x + a}$. **(1)**

Since *(-1)$^2$-4(0)=1>0* (coefficient of the $x^2$ term is *1*; coefficient of the *xy* term is *-1*; and coefficient of the term $y^2$ is *0*); the above degree *2* algebraic equation describes a hyperbola on the *x-y* plane. In L.E. Diskson's book *History of the Theory of Numbers, Vol. II* (see reference **[1]**); there is a wealth of historic information on Diophantine equations of degree *2*; in two or three variables. In reference **[2]**, the reader will find an article published in 2009; and dealing with some special cases of integral points on hyperbolas.

In Section 2, we prove a key proposition, Proposition 1. Proposition 1 lays out the precise (i.e. necessary and sufficient conditions, for a quadratic trinomial with integer coefficients to have two integer roots or zeros. This proposition plays a key role in establishing Theorem 1 in Section 5; the main theorem of this paper. According to Theorem 1, if the integer $a^2 - ab + c$ is not zero (which is to say that the integer *–a* is not a zero of the quadratic trinomial $x^2+bx+c$). Then the graph of the rational function in **(1)** contains exactly *4N* distinct integral points; provided that $a^2$-*ab*+*c* is not equal to an integer or perfect square; or minus an integer square. If on the other hand $a^2$-*ab*+*c*=$k^2$; or *–$k^2$*; where *k* is a positive integer. Then the graph of *y=f(x)* contains exactly *4N-2* distinct integral points. In both cases above, *N* stands for the number of positive integer divisors of $\left|a^2 - ab + c\right|$; divisors not exceeding $\sqrt{\left|a^2 - ab + c\right|}$.

In Section 6, we apply Theorem 1, in the case were $\left|a^2 - ab + c\right| = 1, p, p^2$, or $p_1 p_2$. Where *p* is a prime; and $p_1 p_2$ are distinct primes.

In Section 3, we present an analysis of the graph of $y=f(x)$ in the special case $b^2-4c=0$. In this case the rational function $f(x)$ reduces to $f(x) = \dfrac{(x+d)^2}{x+a}$ ; where $d$ is an integer such that $b=2d$ and $c=d$.

We analyze (in Section 3) this case from the Calculus point of view. As it becomes evident, there are really only two types of graph of $y=f(x)$; in the case $b^2 - 4c = 0$.

See Figures **1, 2, 3, 4,** and **5**

## 2. A Proposition and its Proof

The following proposition, Proposition 1, is the key result used in establishing Theorem 1 in **Section 5**. This proposition and its proof can also be found in reference **[3]** (article authored by this author). But, for reasons of completeness and convenience for the readers of this paper; we include it here. We will make use of a well-known lemma in number theory, (see **[4]**).

**Lemma 1 (Euclid's Lemma)**

Let $n_1$, $n_2$, $n_3$ be nonzero integers such that $n_1$ is a divisor of the product $n_2$, $n_3$. Then if $n_1$ is relatively prime to $n_2$; $n_1$ must divide $n_3$.

**Proposition 1**

*Let a, b, c be integers; with the integer a being nonzero.*

*Consider the quadratic trinomial $g(x) = ax^2 + bx + c$.*

*Then,*

(i) *The trinomial $g(x)$ has either two rational zeros (or roots) or (otherwise) two irrational zeros.*

(ii) *The trinomial $g(x)$ has two rational zeros (or roots) if, and only if, the discriminant $b^2 - 4ac$ is an integer square.*

(iii) *The trinomial $g(x)$ has two integer zeros if, and only if, the following conditions are satisfied:*

$$\begin{cases} b^2 - 4ac = k^2, \text{ for some integer } k. \\ \text{And, the integer } a \text{ is a divisor of both integers } c \text{ and } b \\ (\text{Equivalently, } a \text{ is a divisor of the greatest common divisor of } b \text{ and } c) \end{cases}$$

(iv)  If $a = 1$, then $g(x)$ has two integer zeros if and only if $b^2 - 4c$ is an integer square.

If $a = -1$, then $g(x)$ has two integer zeros if and only if $b^2 + 4c$ is an integer square.

**Proof**

(i)  If $r_1$ and $r_2$ are the zeros of $g(x)$. Then, $r_1 + r_2 = -\dfrac{b}{a}$; (and also $r_1 r_2 = \dfrac{c}{a}$).

Since $-\dfrac{b}{a}$ is a rational number; it is clear from the last equation that if one of $r_1, r_2$ is rational; so is the other.

(ii)  First, suppose that $b^2 - 4ac$ is the square of the square of an integer;

$b^2 - 4ac = D^2$, where $D$ is a nonnegative integer.

The two zeros or roots of the trinomial $g(x)$; are the real numbers $\dfrac{-b + \sqrt{D^2}}{2a}$ and $\dfrac{-b - \sqrt{D^2}}{2a}$;

that is, the numbers $\dfrac{-b + \sqrt{D}}{2a}$ and $\dfrac{-b - \sqrt{D}}{2a}$; which are both rational *since b, D,* and *a* are all integers. Now the converse. Let $r_1$ and $r_2$ be the rational zeros of $g(x)$; and $T$ the discriminant.

We have, $\begin{cases} r_1 = \dfrac{-b + \sqrt{T}}{2a}, r_2 = \dfrac{-b - \sqrt{T}}{2a}; and \\ T = b^2 - 4ac \end{cases}$ **(2)**

We write the rational number $r_1$ in lowest terms:

$\begin{cases} r_1 = \dfrac{u}{v}, where \ u \ and \ v \ are \ relatively \\ prime \ integers \ and \ v \ nonzero \end{cases}$ **(3)**

From (2) and (3) we obtain,

$(2au \vdash bv)^2 = v^2 \cdot T$ **(4)**

If $T = 0$, then $T = b^2$ and we are done: *T* is an integer square.

If T is nonzero, then since *v* is also nonzero; and so by **(4),** so is the integer $2au + bv$. According to **(4),** the nonzero integer square $v^2$ divides the nonzero integer square $(2au + bv)^2$. This implies that positive integer $|v|$ must divide the positive integer $|2au + bv|$ (This last inference can typically be found as an exercise in elementary number theory books. It can also be found in reference **[2]**). Therefore the integer is a divisor of the integer $2au+bv$. And so,

$$\begin{cases} 2au+bv = t \cdot v, \\ \text{for some integer } t \\ \text{or equivalently,} \\ 2au = v \cdot (t-b) \end{cases} \quad (5)$$

If $t=b$, then by **(5)** we get $2au = 0$ which implies $u = 0$ in view of a being nonzero. But then $r_1 = 0$ by **(3)**. Since $g(x)$ has zero as one of its roots. It follows that the constant term $c$ must zero; $c=0$, which implies $T = b^2$ by **(2)**. Again $T$ is an integer square. If $t \neq b$ in **(5)**. Then since all three integers $a, v$, and $t-b$ are nonzero; so must be the integer $u$ by **(5)**. So all four integers $a, u, v, t-b$ are nonzero. By Lemma 1 (Euclid's Lemma) since $v$ is a divisor of the product $2au$; and $v$ is a relatively prime to $u$. It follows that v must be a divisor of the integer $2a$. So that,

$$\begin{cases} 2a = v \cdot w, \text{ for some} \\ \text{nonzero integer } w \end{cases} \quad (6)$$

Combining **(6)** and **(4)** yields,

$$v^2 \cdot (w+b)^2 = v^2 \cdot T ;$$

$$T = (w+b)^2 ,$$

which proves that $T$ is integer square.

The proof is complete. $\square$

*(iii)* Suppose that the trinomial $g(x)$ has two integer roots and zeros.

Then by part (ii) it follows that,

$$\begin{cases} b^2 - 4ac = k^2, \\ \text{for some nonnegative integer } k \end{cases} \quad (7)$$

And the two zeros of $g(x)$ are the numbers

$$\left\{\begin{array}{l} r_1 = \dfrac{-b+k}{2a} \text{ and } \dfrac{-b-k}{2a} = r_2 \\ \text{where } r_1 \text{ and } r_2 \text{ are integers} \end{array}\right\} \quad (8)$$

So, from **(8)** we have,

$$2ar_1 + b = k;$$
$$(2ar_1 + b)^2 = k^2 \quad (9)$$

From **(9)** and **(7)** we get,

$$b^2 + 4abr_1 + 4a^2r_1^2 = b^2 - 4ac; \text{ and since } a \neq 0,$$

$$br_1 + ar_1^2 = -c \quad (10)$$

From **(8)**, we also have,

$$2a(r_1 + r_2) = -2b;$$
$$-a(r_1 + r_2) = b,$$

Which shows that $a$ is a divisor of $b$. Thus,

$$b = \rho \cdot a, \text{ for some integer } \rho \quad (11)$$

By **(11)** and **(10)** we obtain,

$$\rho \cdot a \cdot r_1 + a \cdot r_1^2 = -c;$$

$$-a \cdot r_1 \cdot (\rho + r_1) = c; \text{ which proves}$$

That $a$ is also a divisor of $c$.

We have shown that the integer $a$ is a common divisor of $b$ and $c$.

Now the converse. Suppose that a is a common divisor of both $b$ and $c$; and that $b^2 - 4ac = k^2$, for some $k \in \mathbb{Z}$. Then, $b = a \cdot b_1$ and $c = a \cdot c_1$, for some integers $b_1$ and $c_1$.

And so,

$$a^2 b_1^2 - 4a^2 c_1 = k^2;$$
$$a^2(b_1^2 - 4c_1) = k^2 \quad (12)$$

According to **(12)**; $a^2$ is a divisor of $k^2$ which implies that $a$ is a divisor of $k$. And so $k = ak_1$, for some $k_1 \in \mathbb{Z}$.

Using $b = ab_1, c = ac_1, k = ak_1$; and the formulas in **(8)**. We see that the two zeros $r_1$ and $r_2$ of $g(x)$ are the numbers

$$\left( r_1 = \frac{-b_1 + k_1}{2} \text{ and } r_2 = \frac{-b_1 - k_1}{2} \right) \quad (13)$$

Also from, **(12)**, and $k = ak_1$; we have

$$b_1^2 - 4c_1 = k_1^2 \quad (14)$$

Equation **(14)** shows that $b_1$ and $k_1$ have the same parity; they are either both odd or both even.

Hence by **(13)**, it follows that both rational numbers $r_1$ and $r_2$; are actually integers.

**(iv)** This part is an immediate consequence of part **(iii)**. We omit the details. $\square$

3. **The Case $b^2 - 4c = 0$ : A Calculus Based Analysis**

   Then function $f(x) = \dfrac{x^2 + bx + c}{x + a}$ has domain, $D_f = \{x | x \in \mathbb{R} \text{ and } x \neq -a\} = (-\infty, -a) \cup (-a, +\infty)$.

   Now, suppose that the integers $b$ and $c$ satisfy,

   $$b^2 - 4c = 0;$$
   $$b^2 = 4c$$

   The integer $b$ must be even; $b = 2d$, for some integer $d$.
   Which yields $c = d^2$. We have,

$$f(x) = \frac{x^2 + 2dx + d^2}{x+a};$$

$$f(x) = \frac{(x+d)^2}{x+a}.$$

**Case 1:** $a = d = \dfrac{b}{2}$

In this case the function *f* reduces to,

*f(x)=x+a*

The graph of the $y = f(x)$ is the graph of the straight line with equation $y = x + a$ but with the point $(-a, o)$ removed.

**Figure 1** $y = f(x) = \dfrac{(x+a)^2}{x+a}$; with $a < 0$

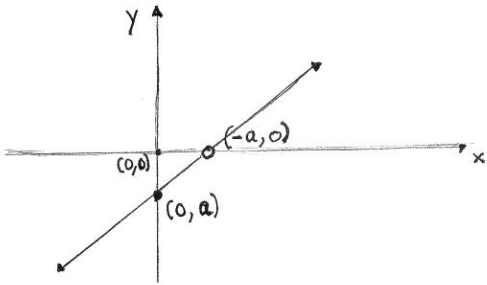

**Figure 2** $y = f(x) = \dfrac{x^2}{x}$

$a = 0$

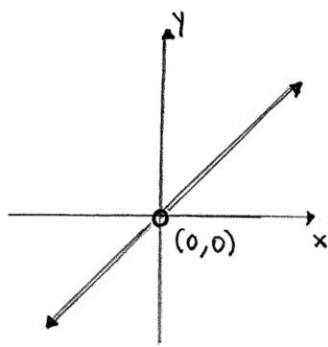

**Figure 3**

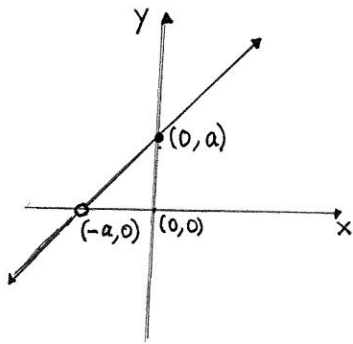

$$y = f(x) = \frac{(x+a)^2}{x+a} \quad \text{with } a < o$$

**Case 2:** $a \neq d = \dfrac{b}{2}$

$$f(x) = \frac{(x+d)^2}{x+a}$$

(i) *Asymptotes*

The graph of $y = f(x)$ has no horizontal asymptote. It has one vertical asymptote, the vertical line $x = -a$. Also, $\lim_{x \to a^-} f(x) = -\infty$, while $\lim_{x \to a^+} = +\infty$.

The graph of $y = f(x)$ also has an oblique asymptote; the slant line with equation $y = x + 2d - a$.

If we perform long division or synthetic division of $(x+d)^2$ with $x+a$; we obtain

$$(x+d)^2 = (x+a)\cdot(x+2d-a) + (d-a)^2;$$

$$f(x) = \frac{(x+d)^2}{x+a} = x + 2d - a + \frac{(d-a)^2}{x+a}; \text{ for } x \neq -a.$$

So that $\lim_{x \to +\infty}[f(x) - (x+2d-a)] = 0 = \lim_{x \to -\infty}[f(x) - (x+2d-a)].$

(ii) *Intercepts*

The point *(-d, 0)* is the *x* intercept on the graph of *y=f(x)*; while the point $\left(0, \dfrac{d^2}{a}\right)$ is the *y intercept*.

(iii) *First Derivatice, Open Intervals of Increase/Decrease, and Points of Relative Extremum.*

We use the quotient rule to calculate the derivative *f'(x)* of *f(x)*:

$$f(x) = \frac{(x+d)^2}{x+a};$$

$$f'(x) = \frac{2(x+d)\cdot(x+a) - (x+d)^2}{(x+a)^2} = \frac{(x+d)\cdot[2(x+a) - (x+d)]}{(x+a)^2}; \qquad (15)$$

$$f'(x) = \frac{(x+d)\cdot[x+2a-d]}{(x+a)^2}.$$

We see from **(15)** that the function *f*, has two critical numbers in its domain $(-\infty, -a)\cup(-a,+\infty)$. These are the numbers *-d* and *–(2a-d)=d-2a*. Note that these two numbers are distinct since *a* and *d* are distinct. Furthermore, when *a<d*, then *–d<-a<d-2a;* and so the function *f* is increasing on the open intervals $(-\infty,-d)$ and $(d-2a,+\infty)$; decreasing on the open intervals *(-d,-a)* and *(-a,d-2a)*.

The point *(-d,0)* on the graph of *y=f(x)*; is a point of relative maximum. While the point *(d-2a, 4(d-a))* is a point of relative maximum. Note that $f(d-2a) = \frac{(2d-2a)^2}{d-a} = 4(d-a).$

On the other hand, if *d<a* . Then $d-2a<-a<-d.$ The function *f* increases on the open intervals $(-\infty, d-2a)$ and $(-d,+\infty)$; and it decreases on the open intervals $(d-2a,-a)$ and $(-a,-d)$. The point $(d-2a, 4(d-a))$ on the graph of $y=f(x)$; is a poiont of relative maximum; while the point $(-d,0)$ is a point of relative minimum when *d<a* and vice-versa when *a>d*.

(iv) *Second derivative, open intervals of concavity, and inflection points*

As we will see below, the graph of $y = f(x)$ has no inflection points. We compute the second derivative. First we write the numerator in **(15)**; in expanded form. $f'(x) = \frac{x^2 + 2ax + d(2a-d)}{(x+a)^2}.$

Applying the quotient rule gives,

$$f''(x) = \frac{2(x+a)\cdot(x+a)^2 - 2(x+a)\cdot[x^2+2ax+d(2a-d)]}{(x+a)^4};$$

$$f''(x) = \frac{2(x+a)^2 - 2[x^2+2ax+d(2a-d)]}{(x+a)^3};$$

$$f''(x) = \frac{2[x^2 + 2ax + a^2 - x^2 - 2ax - d(2a-d)]}{(x+a)^3};$$

$$f''(x) = \frac{2(a-d)^2}{(x+a)^3}; \text{ and since } (a-d)^2 > 0,$$

On account of $a \neq d$. It clear that $f''(x) > o$ on the open interval $(-a, +\infty)$; while $f''(x) < o$ on the open interval $(-\infty, -a)$. Thus, the graph of $y = f(x)$ is concave downwards over the open interval $(-\infty, -a)$; and it concave upwards over the open interval $(-a, +\infty)$. There are no inflection points on the graph of $y = f(x)$.

    (v)    *Two graphs*

**Figure 4**

$a < d$. And so, $-2da+a<-d<-a<d-2a$      $y = f(x) = \dfrac{(x+d)^2}{x+a}$

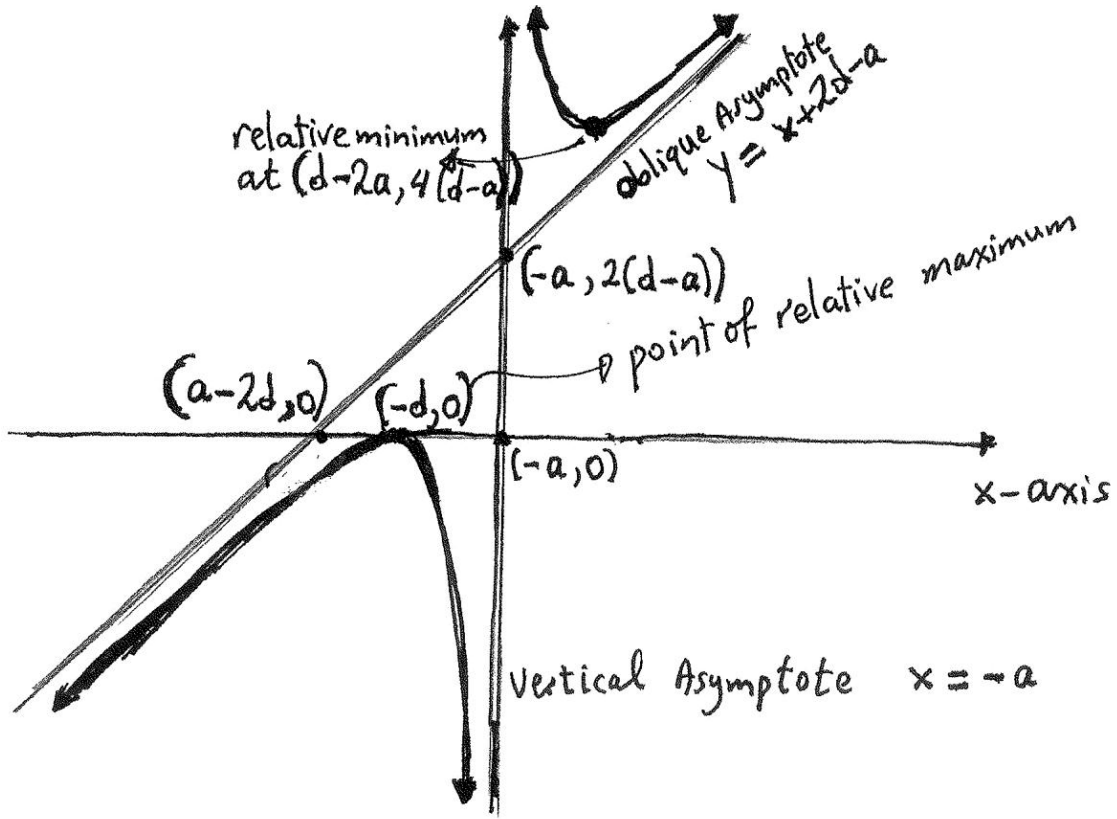

**Figure 5**     $d < a$. And so, $d - 2a < -a < -d < a - 2d$

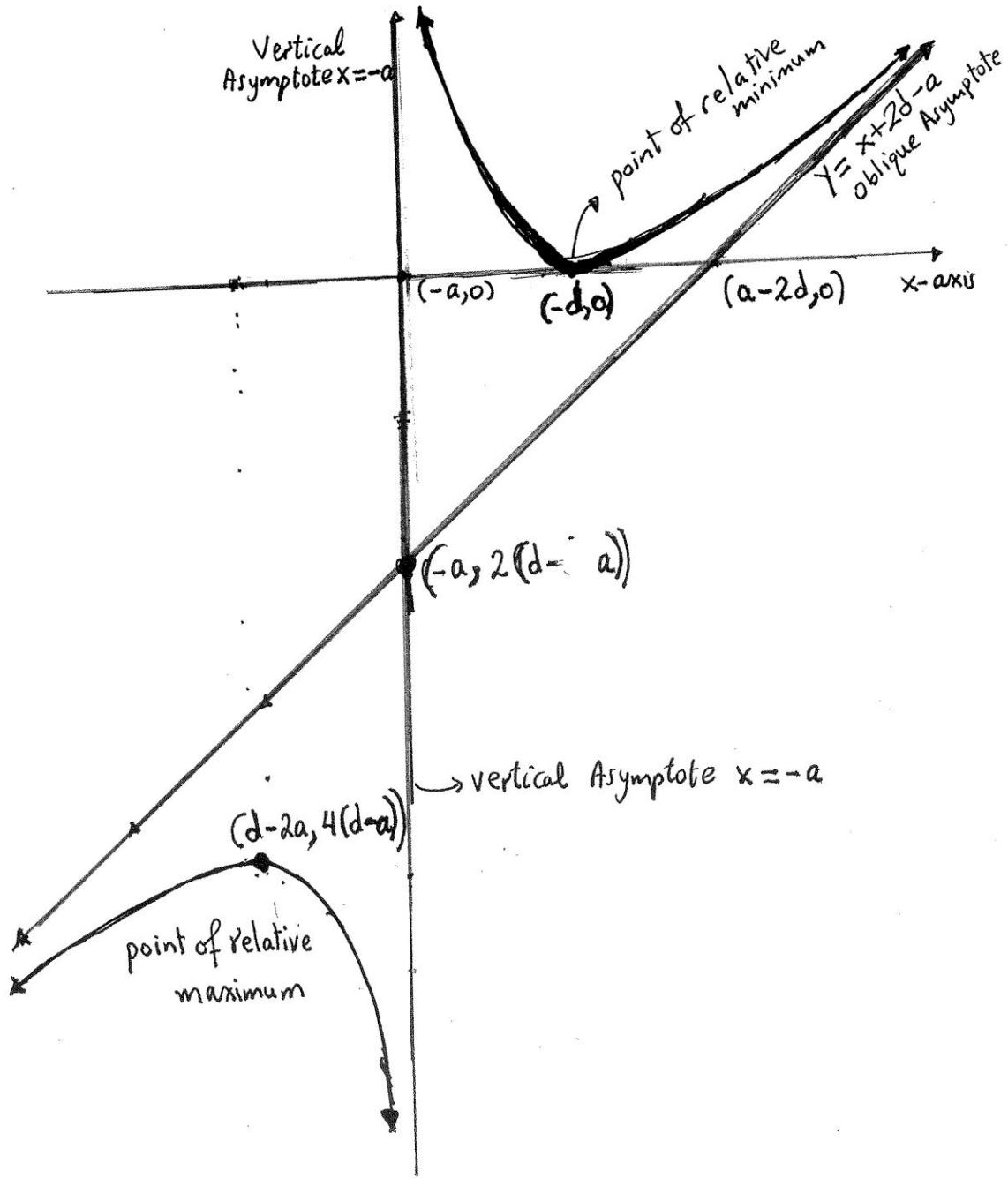

4. **Integral points in the case $b^2 - 4c = o$**

   As we have seen in the previous section; when $b^2 = 4c$, then $b = 2d$ and $c = d^2$; d an integer.

   The function $y = f(x)$ reduces to,

   $$y = f(x) = \frac{(x+d)^2}{x+a}; \text{with } x \neq -a \quad (16)$$

   *Case 1: a=d.* We have $y = x + a$; and so obviously there are infinitely many integral points on the graph of $y = f(x)$. Set of integral points, $S = \{(x, y) | x = t \text{ and } y = t + a; t \in \mathbb{Z}, \text{ and } t \neq -a\}$

   *Case 1: The integers a and d are distinct*

   Equation **(16)** is equivalent to,

   $$\begin{cases} y \cdot (x+a) = (x+d)^2 \\ \text{with } x \neq -a \end{cases} \quad (17)$$

   Clearly, the pair *(-d,0)* is a solution to **(17)**; note that $-d \neq -a$, since $a \neq d$. Now, if $x \neq -d$ in **(17)**. Then $(x+d)^2 > o$. So either *y* and *x+a* are both positive integers; or they are both negative integers.

   So we split the process into two cases.

   *Case 2a:*

   $$\begin{cases} y \cdot (x+a) = (x+d)^2; \\ x, y \text{ are integers such that} \\ y \geq 1, x+a \geq 1 \end{cases} \quad (17a)$$

   Let $\rho$ be the greatest common divisor of *y* and *x+a*.

   Then,

   $$\begin{cases} y = \rho \cdot z_1, \ x + a = \rho \cdot z_2; \\ \text{for some relatively prime positive integers} \\ z_1 \text{ and } z_2 \end{cases} \quad (17b)$$

   From **(17a)** and **(17b)** we obtain

   $$\rho^2 \cdot z_1 \cdot z_2 = (x+d)^2 \quad (17c)$$

According to **(17c)**, the positive integer $\rho^2$ divides the positive integer $(x+d)^2$; which implies that $\rho$ must divide $|x+d|$

$$\begin{cases} |x+d| = \rho \cdot z_3, \\ \text{for some positive integer } z_3 \end{cases} \quad \textbf{(17d)}$$

Combining **(17c)** with **(17d)** yields,

$$z_1 \cdot z_2 = z_3^2 \quad \textbf{(17e)}$$

Since the positive integers $z_1$ and $z_2$ are relatively prime and their product, according to **(17e)**, is equal to a perfect square or integer square; each of them must be an integer square. Recall that more generally; if the product of two relatively prime positive integers is equal to an nth integer power; each of those positive integers must equal an nth integer power. This result can easily be found in number theory books. For examples, see reference **[4].**

Thus by **(17e)** we must have,

$$\begin{cases} z_1 = v_1^2, z_2 = v_2^2, z_3 = v_1 v_2; \\ \text{for relatively prime positive integers} \\ v_1 \text{ and } v_2 \end{cases} \quad \textbf{(17f)}$$

Going back to **(17d)** and **(17b)** to obtain

$$\left( y = \rho \cdot v_1, x = -a + \rho \cdot v_2^2, |x+d| = \rho v_1 v_2 \right) \quad \textbf{(17g)}$$

First, suppose that $a<d$. And so $-a>-d$.

Since $x>-a$. We have $x>-d$, and so by **(17g)** it follows that $|x+d| = x+d = \rho v_1 v_2$. Which further implies by **(17g)** again; that

$$-a + \rho \cdot v_2 v^2 = \rho v_1 v_2 - d;$$
$$\rho \cdot v_2 \cdot (v_1 - v_2) = d - a;$$
$$v_2 \cdot (v_1 - v_2) = \frac{d-a}{\rho}; \text{ with } v_1 - v_2 \geq 1 \text{ and } d - a \geq 1.$$

By setting $m = v_2$ and $v_1 - v_2 = n$; $v_1 = m+n$. We obtain $y = \rho \cdot (m+n)^2$, $x = -a + \rho \cdot m^2$; with m and n being relatively prime integers. Clearly the conditions $y \geq 1$ and $x + a \geq 1$ in **(17a)** are satisfied. Clearly $\gcd(m,n) = 1$; since $\gcd(v_1, v_2) = 1$. Next, suppose that $d < a$; and so $-a < -d$. There are two possibilities to consider in this case. One possibility is $-a < -d < x$. The other possibility is

$-a < x < -d$; which requires that $|-d-(-a)| = |a-d| = a-d \geq 2$. For if $a-d=1$, then there is no integer

in the open interval $(-a, -d)$. Also, clearly $x \neq -d$, as implied by the conditions in **(17a)**.

So, if $-a < -d < x$. Then $x+d > 0; |x+d| = x+d$. So **(17g)** implies $x = -a + \rho \cdot v_2^2 = -d + \rho v_1 v_2$;

$v_2 \cdot (v_2 - v_1) = \dfrac{a-d}{\rho}$; with $v_2 > v_2 - v_1 \geq 1$. By setting $m = v_2$ and $v_2 - v_1 = n$; $v_1 = m - n$. We obtain

$y = \rho \cdot (m-n)^2$, $x = -a + \rho \cdot m^2$; under the conditions $m > n \geq 1$; and with $\gcd(m,n) = 1$. Note that the

conditions $y \geq 1$ and $x+a \geq 1$ in **(17a)** are satisfied. Lastly, if $-a < x < -d$; with $a - d \geq 2$. We have

$x+d < 0$; $x+d \leq -1$ (since $x+d$ is an integer). And then **(17g)** gives,

$$x = -a + \rho \cdot v_2^2 = -d - \rho v_1 v_2$$

And so,

$$v_2 \cdot (v_1 + v_2) = \dfrac{a-d}{\rho}$$

We get, $v_2 = n$, $v_1 = m - n$; with $m$ and $n$ being relatively positive integers with $m > n \geq 1$. We obtain

$y = \rho \cdot (m-n)^2$, $x = -a + \rho \cdot n^2$.

We can now state Result 1, which summarizes the results in Case 2a (the case with $y \geq 1$)

**Result 2**

Consider the rational function $y = f(x) = \dfrac{(x+d)^2}{x+a}$, where a and d are distinct integers.

The set of all integral points (x,y) on the graph of $y = f(x)$; with $y \geq 1$; is a finite set which can be

described as follows:

(i)     If $d > a$. Then the set of all integral points on the graph of $y = f(x)$; with $y \geq 1$, and thus with

$x + a \geq 1$ as well; can be parametrically described in the following manner: $y = \rho \cdot (m+n)^2$,

$x = -a + \rho \cdot m^2$; where m and n are relatively prime positive integers such that $m \cdot n = \dfrac{d-a}{\rho}$;

and $\rho$ is a positive integer divisor of $d - a$.

(ii) If $d < a$. Then set of all integral points on the graph of $y = f(x)$; with $y \geq 1$ can be parametrically described in the following manner. Those integral points with $y \geq 1$ and $x > -d$ are given by, $y = \rho \cdot (m-n)^2$, $x = -a + \rho \cdot m^2$; where m and n are relatively prime integers such that $m > n \geq 1$; and with $m \cdot n = \dfrac{a-d}{\rho}$; $\rho$ a positive divisor of a-d.

And those integral points (such points exist only when $a - d \geq 2$) with $y \geq 1$ and $-a < x < -d$ can be described as follows: $y = \rho \cdot (m-n)^2, x = -a + \rho \cdot n^2$ Where $\rho$ is a positive divisor of a-d; and m, n are relatively prime positive integers such that, $m \cdot n = \dfrac{a-d}{\rho}$; and with $m > n \geq 1$.

**Remark 1**

Note that when $d > a$, and $m = n = 1$, $\rho = d - a$. We obtain the integral point $(x, y) = (d - 2a, \ 4(d-a))$; which is one of the two points of extremum that the curve

$$y = f(x) = \dfrac{(x+d)^2}{x+a}; \text{ as we saw in } \textbf{Section 3}.$$

Now, we go back to **(17)** and we consider the second case.

Case 2b:
$$\begin{bmatrix} y \bullet (x+a) = (x+d)^2; \\ x, y \text{ are integers such that} \\ y \leq -1, \ x + a \leq -1 \end{bmatrix} \quad \textbf{(17h)}$$

If $\rho = \gcd(y, x+a)$. Then,

$$\begin{bmatrix} y = -\rho \cdot z_1 \ , \ x + a = -\rho_2 \cdot z_2; \\ z_1, z_2 \text{ are relatively prime positive integers.} \end{bmatrix} \quad \textbf{(17i)}$$

From **17h** and **17i** one obtains,

$$\begin{bmatrix} z_1 z_2 = z_3^2 \ , \ \rho z_3 = |x+d|; \\ z_3 \text{ a positive integer.} \end{bmatrix} \quad \textbf{(17j)}$$

From this point on the rest of the analysis/proof/procedure for this care (Case2b) is very similar to that of Case2a. We omit the details and state of Result 2.

**Result 2**

Consider the rational function $y = f(x) = \dfrac{(x+d)^2}{x+a}$, where $a$ and $d$ are distinct integers.

The set of all integral points $(x, y)$ on the graph of $y = f(x)$; with $y \leq -1$; is a finite set which can be described as follows:

(i) If $d < a$, then the set of all integral points on the graph of $y = f(x)$; with $y \leq -1$, and thus with $x + a \leq -1$ as well; can be parametrically in the following manner:

$y = -\rho(m+n)^2$, $x = -a - \rho \cdot m^2$; where $\rho$ is a positive divisor of $a - d$; and $m$, $n$ are relatively prime positive integers such that $m \cdot n = \dfrac{a-d}{\rho}$.

(ii) If $a < d$, then the set of all integral points on the graph of $y = f(x)$; with $y \leq -1$; can be parametrically described in the following manner.

Those integral points with $y \leq -1$ and $x < -d$ are given by $y = -\rho(m-n)^2$, $x = -a - \rho m^2$; where $m$ and $n$ are relatively prime positive integers such that $m > n \geq 1$.

And those integral points (such points exist only when $d - a \geq 2$) with $y \leq -1$ and $-d < x < -a$ can be described as follows:

$y = -\rho(m-n)^2$ , $x = -a - \rho n^2$; where $m$ and $n$ are relatively prime positive integers such that $1 \leq n < m$.

We conclude this section by stating the obvious result below.

**Result 3**

Let $f(x) = \dfrac{(x+d)^2}{x+a}$; where $a$ and $d$ are integers.

(i) If $a \neq d$, the only integral point on the graph of $y = f(x)$; with $y-$ coordinate zero; is the point $(-d, 0)$.

(ii) If $a = d$, the set $S$ of integral points on the graph of $y = f(x)$ is the set,

$$S = \{(x, y) \mid x = t,\ y = a + t;\ t \neq -a,\ t \text{ an integer}\}.$$

## 5. The General Case and the Main Theorem

We now go back to $y = \dfrac{x^2 + bx + c}{x + a}$. If $-a$ is a zero of the trinomial $x^2 + bx + c$; i.e. if $a^2 + ab + c = 0$.

Then, the other zero is the integer $-b+a$. So in this case we have

$$y = \dfrac{(x-(a-b))(x+a)}{x+a} = x - (a-b) = x + b - a.$$ All integral points in this case on the above curve are

the points of the form $(x,y)=(t, t+b-a)$; $t$ can be any integer other than $-a$. Next, assume that $a^2 + ab + c$ is not zero. An integral point $(i_1, i_2)$ will be on the curve precisely when,

$$\begin{cases} i_1^2 + i_1 \cdot (b - i_2) + c - a \cdot i_2 = 0, \\ \text{and } i_1 \neq -a \end{cases} \quad \textbf{(18)}$$

According to **(18)**, $i_1$ is an integer zero of the quadratic trinomial (with integer coefficients 1, $b - i_2$, $c - a \cdot i_2$) $t^2 + (b - i_2) + c - a \cdot i_2$; and so the other zero must also be integral. Thus by **Proposition 1(iii)** we must have (and only then):

$$\begin{cases} (b - i_2)^2 - 4 \cdot (c - a \cdot i_2) = K^2, \\ \text{for some nonnegative integer } K. \\ \text{The zeros of (18) are } \dfrac{-(b - i_2) + K}{2} \text{ and } \dfrac{-(b - i_2) - K}{2} \\ i_1 \text{ one of these two zeros.} \end{cases} \quad \textbf{(18a)}$$

From **(18a)**, we have, $i_2^2 - 2 \cdot (b - 2a) \cdot i_2 - (4c + K^2) + b^2 = 0$ **(18b)**

Equation **(18b)** shows that since the integer $i_2$ is one of its two zeros; so must be the other zero. Applying **Proposition 1(iii)** again,

$$\begin{cases} 4(b-2a)^2 - 4\left[b^2 - (4c+K^2)\right] = 4M^2, \\ \text{for some nonnegative integer } M. \\ \text{The zeros of (18c) are the integers} \\ b-2a+M \text{ and } b-2a-M; \\ i_2 \text{ is one of the two integers.} \end{cases} \quad \textbf{(18c)}$$

Taking the equation in **(18c)** further we get,

$$b^2 - 4ab + 4a^2 - b^2 + 4c + K^2 = M^2;$$
$$4(a^2 - ab + c) = M^2 - K^2 = (M-K)(M+K) \quad \textbf{(18d)}$$

Recall that $a^2\text{-}ab\text{+}c$ is nonzero. It is clear from **(18d)** that the nonnegative integers $M$ and $K$ must have the same parity: either they are both odd; or they are both even. Moreover, since $M \geq 0$ and $K \geq 0$; we have $M + K \geq 0$.

Therefore **(18d)** implies $\quad 4 \cdot |a^2 - ab + c| = (M+K) \cdot |M-K|. \quad \textbf{(18e)}$

Since $M, K$ have the same parity, $M+K$ and $M-K$ are both even integers. Also $|a^2 - ab + c|$ is a positive integer, since $a^2\text{-}ab\text{+}c$ is nonzero. Thus we must have,

$$\begin{cases} M + K = 2d_1 \text{ and } |M - K| = 2d_2 \\ \text{where } d_1 \text{ and } d_2 \text{ are positive divisors of } |a^2 - ab + c| \\ \text{such that } d_1 d_2 = |a^2 - ab + c|; \\ \text{and with } 1 \leq d_2 \leq d_1 \end{cases} \quad \textbf{(18f)}$$

Note that $d_2$ cannot exceed $d_1$, since $M + K \geq |M - K|$, on account of the fact that both $M$ and $K$ are nonnegative. The equal sign holds only in the case $K=0$ and $M \geq 1$; or in the case $M = 0$ and $K \geq 1$ (clearly at most one of $K, M$ can be zero, by **(18e)** and the fact that $|a^2 - ab + c| \geq 1$ ). Also, the cases $K = 0,\ M \geq 1$; or $K \geq 1,\ M = 0$ can only occur in the case in which $|a^2 - ab + c|$ is perfect or integer

square. When this happens, $d_2$ can equal $d_1$: $1 \leq d_1 = d_2 = \sqrt{|a^2 - ab + c|}$; and with either $K = 0$ and $M \geq 1$; or alternatively $K \geq 1$ and $M = 0$.

Otherwise, when $|a^2 - ab + c|$ is not a perfect square. Both $M$ and $K$ are positive integers; and so, $M + K > |M - K|$; and $1 \leq d_2 < d_1$. Looking at **(18f)** we see that in all cases, $|a^2 - ab + c| = d_1 d_2 \geq d_2 \cdot d_2 = d_2^2 \geq 1$. And so,

$$1 \leq d_2 \leq \sqrt{|a^2 - ab + c|}$$

Thus all the possible such pairs of divisors of $d_2$ and $d_1$; are obtained by choosing a positive divisor $d_2$ of $|a^2 - ab + c|$; not exceeding $\sqrt{|a^2 - ab + c|}$. Then $d_1 = \dfrac{|a^2 - ab + c|}{d_2}$.

Furthermore from **(18f)** we must have,

$$\begin{cases} M = d_1 + d_2 \text{ and } K = d_1 - d_2; \\ \text{when } 0 \leq K < M \\ \text{While, } M = d_1 - d_2 \text{ and } K = d_1 + d_2; \\ \text{when } 0 \leq M < K \end{cases} \quad \textbf{(18g)}$$

Combining **(18g), (18c),** and **(18d)** we see that, when actually $a^2 - ab + c$ is positive; we must have

$$\begin{cases} M = d_1 + d_2, \ K = d_1 - d_2; \text{ and with} \\ \text{either } i_2 = b - 2a + M = b - 2a + (d_1 + d_2). \\ \text{or alternatively, } i_2 = b - 2a - M = b - 2a - (d_1 + d_2) \end{cases} \quad \textbf{(18h)}$$

And so from **(18h)** and **(18a)**, we also get,

$$\begin{cases} \text{either } i_1 = \dfrac{-[b - (b - 2a) - (d_1 + d_2)] + (d_1 - d_2)}{2}; \\ \text{or } i_1 = \dfrac{-[b - (b - 2a) - (d_1 + d_2)] - (d_1 - d_2)}{2} \\ \text{or } i_1 = \dfrac{-[b - (b - 2a) + (d_1 + d_2)] + (d_1 - d_2)}{2} \\ \text{or } i_1 = \dfrac{-[b - (b - 2a) + (d_1 + d_2)] - (d_1 - d_2)}{2} \end{cases} \quad \textbf{(18i)}$$

The above four possibilities in **(18i)** really reduce to the following four possibilities as a straightforward calculation shows:

$$\begin{cases} i_1 = -a + d_1 \text{ or } i_1 = -a + d_2 \\ \text{or } i_1 = -a - d_2 \text{ or } i_1 = -a - d_1 \end{cases} \quad \textbf{(18j)}$$

Combining **(18j)** with **(18h)**, we see that when $a^2 - ab + c$ is a positive integer; then for each pair of divisors $d_1$ and $d_2$ in **(18f)**. Four integral points are generated. These are,

$$\begin{cases} (i_1, i_2) = (-a + d_1, b - 2a + (d_1 + d_2)), \ (-a + d_2, b - 2a - (d_1 + d_2)), \\ (-a - d_2, b - 2a + (d_1 + d_2)), \ (-a - d_1, b - 2a - (d_1 + d_2)). \end{cases} \quad \textbf{(18k)}$$

By inspection, these four integral points are distinct; with exception of the case where $a^2 - ab + c$ is a perfect square; with the choice $d_1 = d_2 = \sqrt{a^2 - ab + c}$. The four points reduce to two distinct points then.

What happens when $a^2 - ab + c$ is a negative integer? By **(18d)** we must have $0 \leq M < K$. And so by **(18g)** we have $M = d_1 - d_2$ and $K = d_1 + d_2$.

Combining **(18g), (18c),** and **(18d)** one gets either $i_2 = b - 2a + (d_1 - d_2)$ or $i_2 = b - 2a - (d_1 - d_2)$. After that one combines this with **(18a)** to obtain the four integral points corresponding to the choice of divisors $d_1$ and $d_2$ satisfying **(18f)** (see part (iii) of Theorem 1 stated below). We omit the details.

Also, if we choose another pair of divisors $\rho_1$ and $\rho_2$ satisfying **(18f)**; with the pair $(\rho_1, \rho_2)$ being distinct from the pair $(d_1, d_2)$. Then it is clear that the four points generated by the pair $(\rho_1, \rho_2)$; will be distinct from the four points generated by the pair $(d_1, d_2)$. This follows from **(18k)** (and the corresponding formulas in the case $a^2 - ab + c < 0$) and **(18f)**. We omit the details.

We now state Theorem 1.

**Theorem 1**

*Consider the function* $f(x) = \dfrac{x^2 + bx + c}{x + a}$, *with domain all reals except* $-a$; $a, b, c$ *being integers.*

*Then,*

(i) *If* $a^2 - ab + c = 0$. *There are infinitely many integral points on the graph of* $y=f(x)$. *These are the points of the form,* $(x, y) = (t, t + b - a)$, $t$ *an integer,* $t \neq -a$.

(ii) *If* $a^2 - ab + c > 0$. *When there are exactly* $4N$ *integral points on the graph of* $y = f(x)$. *Except in the case where* $a^2 - ab + c$ *is a perfect square, in which case there are exactly 4N-2 such points.*

*Where* $N$ *(in either case) is the number of positive divisors of the integer* $a^2 - ab + c$; *divisors which do not exceed* $\sqrt{a^2 - ab + c}$. *All* $4N$ *integral points can be parametrically described by the formulas,*

$$(x, y) = (-a + d_1, \ b - 2a + (d_1 + d_2)),$$

$$(-a + d_2, \ b - 2a - (d_1 + d_2)),$$

$$(-a - d_2, \ b - 2a + (d_1 + d_2)),$$

$$(-a - d_1, \ b - 2a - (d_1 + d_2)).$$

*Where* $d_1, d_2$ *are positive divisors of* $a^2 - ab + c$; *such that* $1 \leq d_2 \leq d_1$ *and* $d_1 d_2 = a^2 - ab + c$

(iii) *If* $a^2 - ab + c < 0$, *there are exactly* $4N$ *integral points on the graph of* $y = f(x)$, *unless* $a^2 - ab + c$ *is equal to minus a perfect square, in which case there are exactly* $4N - 2$ *such points. Where* $N$ *(in either case) is the number of positive divisors of the integer* $|a^2 - ab + c|$; *divisors which do not exceed* $\sqrt{|a^2 - ab + c|}$. *All* $4N$ *integral points can be described by the formulas,*

$$(x, y) = (-a + d_1, \ b - 2a + (d_1 - d_2)),$$

$$(-a - d_2, \ b - 2a + (d_1 - d_2)),$$

$$(-a+d_2,\ b-2a-(d_1-d_2)),$$

$$(-a-d_1,\ b-2a-(d_1-d_2))$$

Where $d_1$ and $d_2$ are positive integers such that $1 \le d_2 \le d_1$ and $d_1 d_2 = |a^2 - ab + c|$

## 6. An application of Theorem 1: Theorem 2

Theorem 2 below is a direct application of Theorem 1. Consider the cases $|a^2 - ab + c| = 1,\ p,\ \text{or}\ p^2$ where $p$ is a prime.

Since, in Theorem 1, we have $|a^2 - ab + c| = d_1 d_2$, with $1 \le d_2 < d_1$. It follows that when $|a^2 - ab + c| = 1$; then $d_1 = d_2 = 1$

When $|a^2 - ab + c| = p$; then $d_2 = 1$ and $d_1 = p$. When $|a^2 - ab + c| = p^2$; then either $d_2 = 1$ and $d_1 = p$; or $d_1 = d_2 = p$.

We state Theorem 2 without further elaborating.

### Theorem 2

Let $a, b, c$ be integers and consider the function $f(x) = \dfrac{x^2 + bx + c}{x + a}$ with domain $(-\infty, -a)$ U $(-a, +\infty)$.

(a)     If $a^2 - ab + c = 1$, then the graph of $y = f(x)$ contains exactly four integral points. These are
$(-a+1,\ b-2a+2), (-a+1,\ b-2a-2),\ (-a-1,\ b-2a+2),\ \text{and}\ (-a-1,\ b-2a-2)$

(b)     If $a^2 - ab + c = -1$, then the graph of $y = f(x)$ contains exactly two integral points. These are:
$(-a+1,\ b-2a),\ (-a-1,\ b-2a)$

(c)     If $a^2 - ab + c = -p$, $p$ a prime. When the graph of $y = f(x)$ contains exactly four integral points:

$(-a+p, b-2a+(p+1))$, $(-a+1, b-2a-(p+1))$, $(-a-1, b-2a+(p+1))$,

$(-a-p, b-2a-(p+1))$

(d) If $a^2 - ab + c = -p$, $p$ a prime. When the graph of $y = f(x)$ contains exactly four integral points:

$(-a+p, b-2a+(p-1))$, $(-a-1, b-2a+(p-1))$, $(-a+1, b-2a-(p-1))$,

$(-a-p, b-2a-(p-1))$

(e) If $a^2 - ab + c = p^2$, $p$ a prime. When the graph of $y = f(x)$ contains exactly eight integral points:

$(-a+p^2, b-2a+(p^2+1))$, $(-a+1, b-2a-(p^2+1))$,

$(-a-1, b-2a+(p^2+1))$, $(-a-p^2, b-2a-(p^2+1))$,

$(-a+p, b-2a+2p)$, $(-a+p, b-2a-2p)$,

$(-a-p, b-2a+2p)$, $(-a-p, b-2a-2p)$

(f) If $a^2 - ab + c = -p^2$, $p$ a prime. When the graph of $y = f(x)$ contains exactly six integral points. These are:

$(-a+p^2, b-2a+(p^2-1))$, $(-a-1, b-2a+(p^2-1))$,

$(-a+1, b-2a-(p^2-1))$, $(-a-p^2, b-2a-(p^2-1))$,

$(-a+p, b-2a)$, $(-a-p, b-2a)$

**References**

*[1]    L.E. Dickson, History of the Theory of Numbers, Vol II, 803 pp., Chelsea Publishing Company 1992.  ISBN: 0-8128-1935-6.*
*For Pythagorean triangles see pp. 165-190.*
*For Diophantine equations of degree 2, see pp. 341-428.*

*[2]    Konstantine Zelator, "Integral points on hyperbolas: a special case", arXiv:0908.3866, August 2009, 17 pages, no figures.*

*[3]    Konstantine Zelator, "Integer roots of quadratic and cubic polynomials with integer coefficients", arXiv: 1110.6110, October 2011, 14 pages.*

*[4]    W. Sierpinski, Elementary Theory of Numbers, 480 pp. Warsaw, 1964.*
*ISBN: 0-598-52758-3.*
*For Euclid's lemma (Lemma 1) see Theorem 5 on page 14.*
*For Pythagorean triples, see pp. 38-43.*
*Also, see Theorem 8 on page 17.*